\documentclass{amsart}

\usepackage[T1]{fontenc}
\usepackage{lmodern}
\usepackage{microtype}
\usepackage{amsmath,amssymb,amsthm,mathtools}
\IfFileExists{stmaryrd.sty}{\usepackage{stmaryrd}}{}
\usepackage{aliascnt}
\usepackage{enumitem}
\usepackage{xcolor}
\usepackage[colorlinks=true,linkcolor=blue!45!black,citecolor=blue!45!black,urlcolor=blue!45!black]{hyperref}
\usepackage[nameinlink,noabbrev]{cleveref}

\newtheorem{theorem}{Theorem}[section]
\newaliascnt{proposition}{theorem}
\newtheorem{proposition}[proposition]{Proposition}
\aliascntresetthe{proposition}
\newaliascnt{lemma}{theorem}
\newtheorem{lemma}[lemma]{Lemma}
\aliascntresetthe{lemma}
\newaliascnt{question}{theorem}
\newtheorem{question}[question]{Question}
\aliascntresetthe{question}
\newaliascnt{corollary}{theorem}
\newtheorem{corollary}[corollary]{Corollary}
\aliascntresetthe{corollary}
\theoremstyle{definition}
\newaliascnt{definition}{theorem}
\newtheorem{definition}[definition]{Definition}
\aliascntresetthe{definition}
\newaliascnt{construction}{theorem}
\newtheorem{construction}[construction]{Construction}
\aliascntresetthe{construction}
\theoremstyle{remark}
\newaliascnt{remark}{theorem}
\newtheorem{remark}[remark]{Remark}
\aliascntresetthe{remark}

\newcommand{\N}{\mathbb N}
\newcommand{\Ksort}{\mathbb K}
\newcommand{\Bsort}{\mathbb B}
\newcommand{\Prob}{\mathcal P}
\newcommand{\Twisted}{\mathcal Q}
\newcommand{\Two}{\mathbf 2}
\newcommand{\eps}{\varepsilon}
\ifdefined\llbracket
  \newcommand{\angles}[1]{\left\llbracket #1\right\rrbracket}
\else
  \newcommand{\angles}[1]{\mathopen{[\![}#1\mathclose{]\!]}}
\fi

\title{Non-Aware Computable Presentations of Randomizations}
\author{Xuanheng Zhao}
\date{}

\address{School of Mathematics\\
  Nanjing University\\
  Nanjing, Jiangsu 210093, People's Republic of China}

\email{xuanheng21@gmail.com}


\begin{document}

\begin{abstract}
Let $\Two$ be the two-element structure in the empty signature. We build a computable presentation $\Twisted$ of the Borel
randomization $\Two^{[0,1)}$ in which neither constant random variable is a
computable point. Consequently, $\Twisted$ is not aware. This gives a negative
answer to a question of Cuervo Ovalle and Goldbring \cite{CG}. We also construct a countable structure $M$ with no computable presentation but $M^{[0,1)}$ has a computable presentation. This answers another question and establishes the necessity of the awareness condition stated in one of the main theorems in \cite{CG}.
\end{abstract}

\maketitle

\section{Introduction}

Cuervo Ovalle and Goldbring \cite{CG} initiate the effective metric structure
theory of Keisler randomizations \cite{Ke}.

Informally, a Keisler randomization replaces the elements of a first-order
structure by random variables taking values in that structure over a common
atomless probability space. A first-order formula applied to these random
variables determines the event on which it holds, whose probability measures
how often it is satisfied. The randomization records both the random
variables and these events, together with their probabilities. Two random
variables are close when they agree with high probability, and those that
agree almost surely are identified. For the Borel randomizations considered
here, the underlying probability space is the unit interval with Lebesgue
measure. Each element of the original structure is represented by the
corresponding constant random variable, giving a natural copy of the
original structure inside its randomization.

Effective metric structure theory studies how metric structures can be
described and analyzed algorithmically. In a computable presentation, a
countable family of points generates a dense subset on which distances and
basic predicates can be evaluated algorithmically to any prescribed
accuracy. For randomizations, a natural question is whether such a
presentation also provides effective access to the original structure
through its constant random variables.

For a countable structure $M$, Cuervo Ovalle and Goldbring \cite{CG} call a
presentation of its Borel randomization $M^{[0,1)}$ \emph{aware} when the
embedded copy of $M$, viewed as the set of constant random variables, is a
c.e. closed subset relative to the presentation. One of their main results
relates this condition to the decidable presentability of the original
structure.

\begin{theorem}[{\cite[Theorem~3.5]{CG}}]\label{thm:CG}
   A countable structure $M$ has a decidable presentation if and only if $M^{[0,1)}$ has an
aware computable presentation.
\end{theorem}

Intuitively, awareness means that a computable presentation keeps effective
track of the original elements of $M$ inside the space of random variables.
Although its special points may be nonconstant, a single algorithm can use
them to approximate each member of a list consisting of all constant random
variables, to any prescribed accuracy; this is the equivalent criterion in
\cite[Proposition~2.6]{CG}. This explains the role of awareness in
\Cref{thm:CG}: a formula evaluated at constant random variables gives an
event of probability either $0$ or $1$, according to its truth in $M$.
Effective access to the constants therefore allows the truth of formulas in
$M$ to be recovered from the computable data of the randomization.

Their first question about
the additional effectiveness condition for awareness is as follows.

\begin{question}[{\cite[Question~2.11]{CG}}]
For a countable structure $M$, are all computable presentations of $M^{[0,1)}$ aware?
\end{question}

The first purpose of this paper is to answer this question negatively.  The
point is that awareness is defined relative to an externally distinguished
subset of the randomization, namely the constant random variables, and this
subset need not be preserved by an automorphism of the randomization.

For the two-element structure in the empty signature, one may interchange the
two values independently at each sample point. If $A\subseteq[0,1)$ is a
Borel event, interchanging the values precisely on $A$ gives an automorphism
$T_A$ of the randomization. We construct a c.e. open event $A$ such that the
images of the two constant random variables under $T_A$ are not computable in
the standard presentation. We then transport the standard presentation along
$T_A$. Although the transported presentation has the same computable atomic
data as the standard one, the fixed external set of constant random variables
is no longer computably accessible.

\begin{theorem}\label{thm:main}
Let $L$ be the empty signature and let $\Two$ be the $L$-structure with
universe $\{0,1\}$. The Borel randomization $\Two^{[0,1)}$ has a computable
presentation that is not aware. More precisely, there are computably
isomorphic computable presentations $\Prob$ and $\Twisted$ of
$\Two^{[0,1)}$ such that $\Prob$ is the standard induced presentation,
$\Twisted$ is not aware, and neither $\bar 0$ nor $\bar 1$ is a
$\Twisted$-computable point.
\end{theorem}

We also give a second automorphism construction for the countably infinite
pure set.  There every constant remains individually computable, although the
constants cannot be enumerated uniformly as a computable sequence.  These two
examples show that awareness, relative to the fixed set of constant random
variables, is not preserved under computable isomorphisms between
presentations, even when the underlying maps are automorphisms of the
randomization.

There is a second, logically distinct question which is related to \Cref{thm:CG}.

\begin{question}[{\cite[Question~3.6]{CG}}]\label{q2}
    Can one drop the “awareness” assumption in \Cref{thm:CG}?
\end{question}
 
In the final
section we give a construction, based on a computable tree all of
whose paths are computable, which yields the following stronger separation.

\begin{theorem}\label{thm:q36}
There are a computable unary language $L_T$ and a countable $L_T$-structure
$M_T$ such that:
\begin{enumerate}[label=\textup{(\roman*)},leftmargin=2.6em]
  \item $\operatorname{Th}(M_T)$ is decidable and has effective quantifier
        elimination;
  \item $M_T$ has no computable presentation, and hence no decidable
        presentation;
  \item $M_T^{[0,1)}$ has a computable presentation.
\end{enumerate}
Consequently, the answer to \Cref{q2} is negative.
\end{theorem}

The mechanism behind \Cref{thm:q36} is measure-theoretic rather than
automorphic.  A noncomputable point of the probability space can be used to
produce, pointwise, a listing that cannot be produced computably.  The
computable points form a null set, so modifying the random variables there
makes them everywhere well-defined without changing any event in the
probability algebra.  Effective quantifier elimination then turns the
resulting almost-everywhere listing into a computable dense presentation of
the randomization.

We use standard notation from continuous logic, see \cite{Ben}. For effective metric structure theory and basic concepts in computability, see \cite{Fr}.

The paper is organized as follows.  In \Cref{sec:preliminaries} we recall the
presentation conventions for randomizations.  In \Cref{sec:first-construction}
we prove \Cref{thm:main} by a direct diagonalization and an automorphic change
of presentation.  The countably infinite pure-set variant is given in
\Cref{sec:infinite-pure-set}.  Finally, \Cref{sec:q36} proves
\Cref{thm:q36}.

\section{Preliminaries}\label{sec:preliminaries}

We recall the definitions needed for the construction. Our terminology and
presentation conventions follow \cite[Section~2]{CG}.

Fix a countable first-order language $L$ and a countable $L$-structure $M$.
The Borel randomization $M^{[0,1)}$ is a two-sorted continuous structure.

\begin{definition}\label{def:randomization}
The \emph{random-variable sort} $\Ksort$ of the \emph{randomization pre-structure} consists
of the Borel measurable functions $f\colon[0,1)\to M$, where $M$ carries the
discrete sigma-algebra. The \emph{event sort} $\Bsort$ is the probability algebra of
$([0,1),\mu)$, where $\mu$ is Lebesgue measure. For each $L$-formula
$\varphi(x_1,\ldots,x_n)$, the randomization language has a function symbol
\[
 \angles{\varphi}\colon\Ksort^n\longrightarrow\Bsort
\]
interpreted by
\[
 \angles{\varphi(f_1,\ldots,f_n)}
 =\{\omega\in[0,1):
 M\models\varphi(f_1(\omega),\ldots,f_n(\omega))\},
\]
where the set on the right is identified with its equivalence class modulo
null sets. The pseudometrics are
\[
 d_{\Ksort}(f,g)=\mu\angles{f\neq g}
 \quad\text{and}\quad
 d_{\Bsort}(E,F)=\mu(E\mathbin\triangle F).
\]
One then passes to the metric quotients and completions. We continue to write
$M^{[0,1)}$, $\Ksort$, and $\Bsort$ for the completed structure and its two
sorts.
\end{definition}

Each $a\in M$ is identified with the corresponding constant random variable
$\bar a(\omega)=a$. Thus $M$ is viewed as a subset of $\Ksort$. Distinct
constant random variables have distance $1$, so this embedded copy of $M$ is a
closed discrete subset of $\Ksort$.

When $M=\Two$ and $A\subseteq[0,1)$ is Borel, let
\[
 \chi_A(\omega)=
 \begin{cases}
  1,&\omega\in A,\\
  0,&\omega\notin A.
 \end{cases}
\]
Changing $A$ on a null set does not change the corresponding point of
$\Ksort$, so we also write $\chi_A$ when $A$ denotes an event in the
probability algebra.

\begin{definition}\label{def:presentation}
Let $N$ be a separable metric structure. A \emph{presentation}
$N^\#=(N,(a_i)_{i\in\N})$ is a countable sequence of \emph{special
points} whose closure under the function symbols is dense in $N$. Points
obtained from the special points by applying the function symbols are called
\emph{generated points}.
In the many-sorted setting, one specifies such a sequence in each sort.
\end{definition}

There are no function symbols with codomain $\Ksort$ in the randomization
language. Consequently, in the $\Ksort$-sort, the special points and generated
points coincide.

\begin{definition}\label{def:computable-presentation}
A presentation $N^\#$ in a computable signature is \emph{computable} if there is an
algorithm which, given an $n$-ary predicate symbol $P$, a tuple $\bar a$ of
generated points, and a positive rational $\eps$, returns a rational $r$ such
that
\[
 |P^N(\bar a)-r|<\eps.
\]
In particular, the metric on generated points is uniformly computable.
\end{definition}

\begin{definition}\label{def:computable-data}
Let $N^\#$ be a presentation of $N$.
\begin{enumerate}[label=\textup{(\arabic*)},leftmargin=2.6em]
\item A point $x\in N$ is \emph{$N^\#$-computable} if there is an
  algorithm which,
  given a positive rational $\eps$, returns a generated point $q$ such that
  $d(x,q)<\eps$.
\item A sequence $(x_n)_{n\in\N}$ is an \emph{$N^\#$-computable
  sequence} if the
  approximation algorithm in item~\textup{(1)} is uniform in $n$.
\item If $N_0^\#$ and $N_1^\#$ are presentations and $F\colon N_0\to N_1$ is
  an embedding, then $F$ is a \emph{computable map} from $N_0^\#$ to
  $N_1^\#$ if,
  given a generated point $p$ of $N_0^\#$ and a positive rational $\eps$, one
  can compute a generated point $q$ of $N_1^\#$ with $d(F(p),q)<\eps$.
  A bijective $F$ is a \emph{computable isomorphism} if both $F$ and
  $F^{-1}$ are
  computable maps.
\end{enumerate}
\end{definition}

Fix an effective enumeration $(p_i)_{i\in\N}$ of the $\Ksort$-special points
of a computable presentation $\Prob$.

\begin{definition}\label{def:fast-name}
A partial computable function $\Phi\colon\N\rightharpoonup\N$ is a
\emph{fast name} for $f\in\Ksort$ if $\Phi$ is total and
\[
 d_{\Ksort}(p_{\Phi(n)},f)<2^{-n}
 \qquad\text{for every }n\in\N.
\]
\end{definition}

\begin{lemma}\label{lem:fast-name}
A point $f\in\Ksort$ is $\Prob$-computable if and only if it has a fast name.
\end{lemma}

\begin{proof}
A fast name directly gives the approximations required in
\Cref{def:computable-data}. Conversely, if $f$ is $\Prob$-computable, apply
its approximation procedure with error $2^{-n}$ and return the index of the
resulting $\Ksort$-generated point. Since generated and special points coincide
in $\Ksort$, this gives a total computable fast name.
\end{proof}

\begin{definition}\label{def:ce-closed}
Let $N^\#$ be a presentation of $N$. A closed set $X\subseteq N$ is \emph{c.e.
closed} in $N^\#$ if the set of rational open balls
\[
 B(p;r)=\{x\in N:d(x,p)<r\},
\]
where $p$ is generated and $r>0$ is rational, that intersect $X$ is
computably enumerable.
\end{definition}

\begin{definition}\label{def:aware}
A presentation of $M^{[0,1)}$ is \emph{aware} if the embedded set of
constant random variables
\[
 \{\bar a:a\in M\}\subseteq\Ksort
\]
is c.e. closed in that presentation.
\end{definition}

\begin{proposition}[{\cite[Proposition~2.6]{CG}}]\label{prop:aware-criterion}
Let $\mathcal S$ be a computable presentation of $M^{[0,1)}$. Then
$\mathcal S$ is aware if and only if there is an enumeration
$(a_n)_{n\in\N}$ of all constant random variables that is an
$\mathcal S$-computable sequence.
\end{proposition}

Let $\Two$ be the two-element structure in the empty signature. Its standard
presentation is decidable. Let
\[
 \Prob=\bigl(\Two^{[0,1)}\bigr)^\#
\]
be the induced presentation from \cite[Section~2.3]{CG}. Its
$\Ksort$-special points are binary-valued simple functions whose level sets
belong to the Boolean algebra generated by rational intervals, and its
$\Bsort$-special points form the standard presentation of the probability
algebra. In particular, $\bar0$ and $\bar1$ occur effectively among the
$\Ksort$-special points.

\begin{lemma}\label{lem:P-computable}
The presentation $\Prob$ is computable and aware. Moreover, uniformly in a
$\Ksort$-special point $p$, a rational interval $J$, and a positive rational
$\eps$, one can compute a rational $r$ such that
\[
 \left|r-\mu\bigl(J\cap\angles{p=\bar1}\bigr)\right|<\eps.
\]
\end{lemma}

\begin{proof}
Computability is a special case of the fact that the induced randomization of
a decidable presentation is computable \cite[Corollary~3.4]{CG}. It can also
be seen directly. Given finitely many binary simple functions, refine their
finite rational partitions to a common finite partition. On each cell, every
input has a fixed value in $\{0,1\}$, and the truth of any formula of the empty
signature in the finite structure $\Two$ is decidable. Hence every atomic
event of the randomization is a finite Boolean combination of rational
intervals with computable measure. The same argument, applied to
$J\cap\angles{p=\bar1}$, gives the asserted uniform computation.

For awareness, let $p$ be a binary simple function and let $r>0$ be rational.
Since
\[
 d_{\Ksort}(p,\bar0)=\mu\angles{p=\bar1}
 \quad\text{and}\quad
 d_{\Ksort}(p,\bar1)=\mu\angles{p=\bar0},
\]
the ball $B(p;r)$ meets $\{\bar0,\bar1\}$ if and only if
\[
 \mu\angles{p=\bar0}>1-r
 \quad\text{or}\quad
 \mu\angles{p=\bar1}>1-r.
\]
Both measures are computable, and a strict inequality between a computable
real and a rational is c.e. Thus the rational open balls meeting the constants
can be enumerated.
\end{proof}

\section{The proof of Theorem \ref{thm:main}}\label{sec:first-construction}

We first construct a c.e. open event $A$ such that neither $\chi_A$ nor
$\chi_{A^c}$ is a $\Prob$-computable point.

Fix an effective enumeration $(\Phi_e)_{e\in\N}$ of all partial computable
functions from $\N$ to $\N$. The requirement $\mathcal R_e$ is
\[
 \mathcal R_e:\qquad
 \text{if $\Phi_e$ is a fast name for some $f_e$, then
 $f_e\neq\chi_A$ and $f_e\neq\chi_{A^c}$.}
\]

For $m\in\N$, let
\[
 I_m=(2^{-(m+1)},2^{-m}).
\]
These are pairwise disjoint rational intervals. Assign to $\mathcal R_e$ the
two private intervals
\[
 J_e^0=I_{2e},\qquad J_e^1=I_{2e+1},
\]
and set $$\delta_e^k=\mu(J_e^k)>0.$$
Effectively choose $N_e$ so large that,
with $\eps_e=2^{-N_e}$,
\begin{equation}\label{eq:epsilon-choice}
 2\eps_e<\frac14\min\{\delta_e^0,\delta_e^1\}.
\end{equation}

\begin{construction}\label{cons:A}
At stage $s$, for every $e\le s$ that has not yet acted, simulate
$\Phi_e(N_e)$ for $s$ steps. If
\[
 \Phi_e(N_e)\downarrow=i,
\]
use \Cref{lem:P-computable} to compute rationals $r_e^0,r_e^1$ such that
\begin{equation}\label{eq:r-approx}
 \left|r_e^k-
 \mu\bigl(J_e^k\cap\angles{p_i=\bar1}\bigr)\right|<\eps_e
 \qquad(k=0,1).
\end{equation}
Then make the following permanent decisions.

On $J_e^0$, define $A$ to have the value opposite to the approximate majority
value of $p_i$:
\begin{equation}\label{eq:J0-action}
 A\cap J_e^0=
 \begin{cases}
  \varnothing,&r_e^0\ge\delta_e^0/2,\\
  J_e^0,&r_e^0<\delta_e^0/2.
 \end{cases}
\end{equation}
On $J_e^1$, define $A$ to have the same value as the approximate majority
value of $p_i$:
\begin{equation}\label{eq:J1-action}
 A\cap J_e^1=
 \begin{cases}
  J_e^1,&r_e^1\ge\delta_e^1/2,\\
  \varnothing,&r_e^1<\delta_e^1/2.
 \end{cases}
\end{equation}
Finally, declare that $\mathcal R_e$ has acted. Whenever an interval is to be
included in $A$, enumerate that interval at this stage. If $\mathcal R_e$
never acts, neither of its private intervals is included in $A$.
\end{construction}

Since each private interval is controlled by a unique requirement, the
decisions are compatible. The resulting set $A$ is a c.e. union of rational
open intervals and hence is a c.e. open, in particular Borel, subset of
$[0,1)$.

\begin{proposition}\label{prop:diagonalization}
For the event $A$ constructed above, neither $\chi_A$ nor $\chi_{A^c}$ is a
$\Prob$-computable point.
\end{proposition}

\begin{proof}
Fix $e$ and suppose that $\Phi_e$ is a fast name for some $f_e\in\Ksort$.
Then $\Phi_e(N_e)$ converges, so $\mathcal R_e$ acts. Write
\[
 i=\Phi_e(N_e),\qquad p=p_i,
\]
and, for $k=0,1$, set
\[
 \alpha_e^k=
 \mu\bigl(J_e^k\cap\angles{f_e=\bar1}\bigr).
\]
Because $\Phi_e$ is a fast name,
\[
 d_{\Ksort}(p,f_e)<\eps_e.
\]
For binary-valued random variables, the map
$f\mapsto\angles{f=\bar1}$ is $1$-Lipschitz from $\Ksort$ to $\Bsort$.
Consequently,
\[
 \begin{split}
 &\left|
 \mu\bigl(J_e^k\cap\angles{p=\bar1}\bigr)
 -\mu\bigl(J_e^k\cap\angles{f_e=\bar1}\bigr)
 \right|\\
 &\hspace{4cm}\le
 d_{\Bsort}\bigl(\angles{p=\bar1},\angles{f_e=\bar1}\bigr)
 \le d_{\Ksort}(p,f_e)<\eps_e.
 \end{split}
\]
Combining this with \Cref{eq:r-approx}, we obtain
\begin{equation}\label{eq:alpha-r}
 |\alpha_e^k-r_e^k|<2\eps_e
 \qquad(k=0,1).
\end{equation}

We first diagonalize against $\chi_A$ on $J_e^0$. If
$r_e^0\ge\delta_e^0/2$, then $A\cap J_e^0=\varnothing$ by
\Cref{eq:J0-action}, and hence
\[
 \mu\bigl(J_e^0\cap\angles{f_e\neq\chi_A}\bigr)
 =\alpha_e^0
 >\frac{\delta_e^0}{2}-2\eps_e
 >\frac{\delta_e^0}{4}.
\]
If $r_e^0<\delta_e^0/2$, then $A\cap J_e^0=J_e^0$, and hence
\[
 \mu\bigl(J_e^0\cap\angles{f_e\neq\chi_A}\bigr)
 =\delta_e^0-\alpha_e^0
 >\delta_e^0-r_e^0-2\eps_e
 >\frac{\delta_e^0}{4}.
\]
The final inequalities in both cases use \Cref{eq:epsilon-choice}. Therefore
\begin{equation}\label{eq:not-A}
 d_{\Ksort}(f_e,\chi_A)>\frac{\delta_e^0}{4}>0.
\end{equation}

On $J_e^1$, the definition of $A$ makes $A^c$ have the value opposite to the
approximate majority value of $p$. If $r_e^1\ge\delta_e^1/2$, then
$A^c\cap J_e^1=\varnothing$, so
\[
 \mu\bigl(J_e^1\cap\angles{f_e\neq\chi_{A^c}}\bigr)
 =\alpha_e^1
 >\frac{\delta_e^1}{4}.
\]
If $r_e^1<\delta_e^1/2$, then $A^c\cap J_e^1=J_e^1$, so
\[
 \mu\bigl(J_e^1\cap\angles{f_e\neq\chi_{A^c}}\bigr)
 =\delta_e^1-\alpha_e^1
 >\frac{\delta_e^1}{4}.
\]
Consequently,
\begin{equation}\label{eq:not-Ac}
 d_{\Ksort}(f_e,\chi_{A^c})>\frac{\delta_e^1}{4}>0.
\end{equation}
This verifies $\mathcal R_e$.

By \Cref{lem:fast-name}, every $\Prob$-computable point has a fast name
$\Phi_e$ for some $e$. Equations \eqref{eq:not-A} and \eqref{eq:not-Ac}
therefore show that neither $\chi_A$ nor $\chi_{A^c}$ is
$\Prob$-computable.
\end{proof}

\begin{remark}\label{rem:simpler-choice}
For the counterexample alone, it would suffice to choose a Borel event $A$ for
which $\chi_A$ is not $\Prob$-computable. Indeed, the global interchange of
$0$ and $1$ is a computable automorphism of the standard presentation and
sends $\chi_A$ to $\chi_{A^c}$, so one of these points is computable exactly
when the other is. The construction above is retained because it gives a
direct diagonalization against all computable random variables and
makes both noncomputability statements explicit.
\end{remark}

Let $\tau$ be the nontrivial automorphism of $\Two$, so that
$\tau(0)=1$ and $\tau(1)=0$. For $\omega\in[0,1)$, define
\[
 \tau_\omega=
 \begin{cases}
  \tau,&\omega\in A,\\
  \operatorname{id}_{\Two},&\omega\notin A.
 \end{cases}
\]
For a Borel random variable $f$ in the randomization pre-structure, set
\begin{equation}\label{eq:TA}
 T_A(f)(\omega)=\tau_\omega(f(\omega)),
\end{equation}
and let $T_A$ act as the identity on $\Bsort$.

\begin{lemma}\label{lem:TA-automorphism}
The map $T_A$ descends to an involutive automorphism of
$\Two^{[0,1)}$. For every $L$-formula
$\varphi(x_1,\ldots,x_n)$ and every $f_1,\ldots,f_n\in\Ksort$,
\begin{equation}\label{eq:event-invariance}
 \angles{\varphi(T_Af_1,\ldots,T_Af_n)}
 =\angles{\varphi(f_1,\ldots,f_n)}.
\end{equation}
Furthermore,
\begin{equation}\label{eq:constant-images}
 T_A(\bar0)=\chi_A,
 \qquad
 T_A(\bar1)=\chi_{A^c}.
\end{equation}
\end{lemma}

\begin{proof}
Since $A$ is Borel, $T_A(f)$ is Borel whenever $f$ is Borel. For every
$\omega$, the map $\tau_\omega$ is an automorphism of the empty-signature
structure $\Two$. Hence, for Borel random variables and pointwise in
$\omega$,
\[
 \Two\models\varphi(f_1(\omega),\ldots,f_n(\omega))
\]
if and only if
\[
 \Two\models\varphi(\tau_\omega f_1(\omega),\ldots,
 \tau_\omega f_n(\omega)).
\]
This proves \Cref{eq:event-invariance} on the pre-structure. It also gives
\[
 f(\omega)=g(\omega)
 \quad\Longleftrightarrow\quad
 T_Af(\omega)=T_Ag(\omega),
\]
so $T_A$ preserves the $\Ksort$-pseudometric. It is the identity on
$\Bsort$, preserves every distinguished formula-event function, and satisfies
$T_A^2=\operatorname{id}$. It therefore descends to the metric quotient and
extends uniquely to an involutive automorphism of the completion. Equality
\eqref{eq:event-invariance} on the completion follows by continuity.
Finally, \Cref{eq:constant-images} is immediate from the definition.
\end{proof}

Let $(p_i)_{i\in\N}$ be the $\Ksort$-special points of $\Prob$, and let
$(b_j)_{j\in\N}$ be its $\Bsort$-special points. Define a new presentation
$\Twisted$ by taking
\[
 q_i=T_A(p_i)\qquad(i\in\N)
\]
as its $\Ksort$-special points and retaining $(b_j)_{j\in\N}$ as its
$\Bsort$-special points. Since $T_A$ is an automorphism, $(q_i)$ is dense in
$\Ksort$, so $\Twisted$ is a presentation of the same randomization.

\begin{lemma}\label{lem:Q-computable}
The presentation $\Twisted$ is computable. Moreover,
$T_A\colon\Prob\to\Twisted$ is a computable isomorphism.
\end{lemma}

\begin{proof}
For all $i,j$,
\[
 d_{\Ksort}(q_i,q_j)
 =d_{\Ksort}(T_Ap_i,T_Ap_j)
 =d_{\Ksort}(p_i,p_j).
\]
For every $L$-formula $\varphi$ and every tuple of indices $\bar i$,
\Cref{eq:event-invariance} gives
\[
 \angles{\varphi(q_{\bar i})}
 =\angles{\varphi(p_{\bar i})}.
\]
The $\Bsort$-special points are unchanged, and $T_A$ is the identity on
$\Bsort$. By induction on terms, every generated $\Bsort$-term in
$\Twisted$ denotes the same event as the term with the same syntactic code in
$\Prob$. Thus there is an effective code-preserving translation from generated
points of $\Twisted$ to generated points of $\Prob$ under which all metric and
predicate values agree. The algorithms computing the atomic data of $\Prob$
therefore compute the atomic data of $\Twisted$.

On the $\Ksort$, $T_A$ sends the special point $p_i$ to the special point
$q_i$ with the same index and error $0$. On the $\Bsort$, it sends every
generated term to the term with the same syntactic code. Hence $T_A$ is a
computable map from $\Prob$ to $\Twisted$. Since $T_A$ is an involution, the
same index-preserving procedures compute its inverse. Therefore $T_A$ is a
computable isomorphism.
\end{proof}

\begin{remark}\label{rem:legitimacy}
The event $A$ is deliberately noncomputable as a point of the standard
presentation, but this does not prevent $\Twisted$ from being a legitimate
computable presentation. A presentation is an indexed dense family of actual
elements of the metric structure; its special points are not required to be
computable in a previously fixed presentation. Computability is determined by
the metric and predicate values on the newly indexed generated points, and
those values are computable here because $T_A$ preserves them exactly.
The map $T_A$ is computable from $\Prob$ to $\Twisted$, but not from
$\Prob$ to itself: it sends the $\Prob$-computable point $\bar0$ to the
noncomputable point $\chi_A$.
\end{remark}

\begin{proof}[Proof of \Cref{thm:main}]
The presentation $\Twisted$ is computable by \Cref{lem:Q-computable}. Suppose,
toward a contradiction, that $\Twisted$ is aware. By
\Cref{prop:aware-criterion}, there is a $\Twisted$-computable sequence
$(c_n)_{n\in\N}$ whose range is exactly $\{\bar0,\bar1\}$. Since
$T_A^{-1}\colon\Twisted\to\Prob$ is a computable isometry, the sequence
\[
 f_n=T_A^{-1}(c_n)
\]
is $\Prob$-computable, and its range is
$\{\chi_A,\chi_{A^c}\}$ by \Cref{eq:constant-images}.

In particular, $f_0$ is a $\Prob$-computable point belonging to
$\{\chi_A,\chi_{A^c}\}$, contrary to \Cref{prop:diagonalization}.
Hence $\Twisted$ is not aware.

Finally, if $\bar0$ were a $\Twisted$-computable point, then applying the
computable inverse $T_A^{-1}$ would make $\chi_A$ a $\Prob$-computable point.
Similarly, computability of $\bar1$ in $\Twisted$ would imply computability of
$\chi_{A^c}$ in $\Prob$. Thus neither constant is
$\Twisted$-computable.
\end{proof}

\begin{corollary}\label{cor:not-invariant}
Awareness, relative to the fixed embedded copy of the original structure, is
not preserved under computable isomorphisms between presentations of the
same randomization.
\end{corollary}

\begin{proof}
The standard presentation $\Prob$ is aware by \Cref{lem:P-computable}, while
$\Twisted$ is not aware, and \Cref{lem:Q-computable} gives a computable
isomorphism $T_A\colon\Prob\to\Twisted$ whose underlying map is an
automorphism of the randomization.
\end{proof}

\section{A second construction: the infinite pure set}
\label{sec:infinite-pure-set}

The preceding construction has the stronger feature that neither constant
random variable is computable in the new presentation. We now give a shorter
counterexample using a different underlying structure. In this example every
constant remains individually computable, but the constants cannot be
enumerated uniformly as a computable sequence.

Let $M_\infty$ be the structure with universe $\N$ in the empty signature.
Thus every permutation of $\N$ is an automorphism of $M_\infty$.

\begin{lemma}\label{lem:infinite-involution}
There is an involution $q\colon\N\to\N$ such that its fixed-point set
\[
 F=\{n\in\N:q(n)=n\}
\]
is not c.e.
\end{lemma}

\begin{proof}
Fix a c.e. noncomputable set $K\subseteq\N$. On each pair
$\{2e,2e+1\}$, let $q$ interchange the two elements if $e\in K$ and fix
both elements if $e\notin K$. Explicitly,
\[
 q(2e)=
 \begin{cases}
  2e+1,&e\in K,\\
  2e,&e\notin K,
 \end{cases}
 \qquad
 q(2e+1)=
 \begin{cases}
  2e,&e\in K,\\
  2e+1,&e\notin K.
 \end{cases}
\]
Then $q^2=\operatorname{id}_{\N}$ and
\[
 F=\{2e,2e+1:e\notin K\}.
\]
If $F$ were c.e., then
$\N\setminus K=\{e:2e\in F\}$ would be c.e. Since $K$ is c.e., this would
make $K$ computable, a contradiction.
\end{proof}

Put
\[
 I_0=[0,\tfrac12),
 \qquad
 I_1=[\tfrac12,1).
\]
For a Borel random variable $f\colon[0,1)\to\N$, define
\begin{equation}\label{eq:Theta-infinite}
 \Theta(f)(\omega)=
 \begin{cases}
  f(\omega),&\omega\in I_0,\\
  q(f(\omega)),&\omega\in I_1,
 \end{cases}
\end{equation}
and let $\Theta$ act as the identity on $\Bsort$.

\begin{lemma}\label{lem:Theta-infinite-automorphism}
The map $\Theta$ descends to an involutive automorphism of
$M_\infty^{[0,1)}$. Moreover, for every formula
$\varphi(x_1,\ldots,x_n)$ in the empty signature and all
$f_1,\ldots,f_n\in\Ksort$,
\begin{equation}\label{eq:Theta-infinite-invariance}
 \angles{\varphi(\Theta f_1,\ldots,\Theta f_n)}
 =\angles{\varphi(f_1,\ldots,f_n)}.
\end{equation}
\end{lemma}

\begin{proof}
The function $\Theta(f)$ is Borel because $\N$ has the discrete
sigma-algebra. At each $\omega$, the fiber map used in
\Cref{eq:Theta-infinite} is either $\operatorname{id}_{\N}$ or $q$, hence an
automorphism of $M_\infty$. Formula satisfaction is therefore preserved
pointwise, which gives \Cref{eq:Theta-infinite-invariance}. Applying this to
equality shows that $\Theta$ preserves the $\Ksort$-pseudometric. It is the
identity on $\Bsort$ and satisfies $\Theta^2=\operatorname{id}$, so it
descends to the metric quotient and extends uniquely to an involutive
automorphism of the completion.
\end{proof}

Let
\[
 \Prob_\infty=\bigl(M_\infty^{[0,1)}\bigr)^\#
\]
be the standard induced presentation. The standard presentation of
$M_\infty$ is decidable, so $\Prob_\infty$ is computable by
\cite[Corollary~3.4]{CG}. Let $(r_i)_{i\in\N}$ be its $\Ksort$-special
points and $(s_j)_{j\in\N}$ its $\Bsort$-special points. Since the induced
presentation effectively lists all rational simple random variables, there
is a computable function $c\colon\N\to\N$ such that
$$r_{c(n)}=\bar n
 \qquad(n\in\N).$$

Define a presentation $\Twisted_\infty$ by taking
\[
 u_i=\Theta(r_i)\qquad(i\in\N)
\]
as its $\Ksort$-special points and retaining $(s_j)_{j\in\N}$ as its
$\Bsort$-special points.

\begin{lemma}\label{lem:Q-infinite-computable}
The presentation $\Twisted_\infty$ is computable, and
$\Theta\colon\Prob_\infty\to\Twisted_\infty$ is a computable isomorphism.
\end{lemma}

\begin{proof}
By \Cref{lem:Theta-infinite-automorphism}, $(u_i)_{i\in\N}$ is dense in
$\Ksort$. Equation \eqref{eq:Theta-infinite-invariance} shows that every
formula-event evaluated on a tuple of the $u_i$ is exactly the same event as
the corresponding formula-event evaluated on the tuple of the $r_i$. The
$\Bsort$-special points are unchanged. Consequently, generated points with
the same syntactic codes in $\Prob_\infty$ and $\Twisted_\infty$ have the
same metric and predicate values. The algorithms for $\Prob_\infty$ therefore
compute the atomic data of $\Twisted_\infty$.

The map $\Theta$ sends each $r_i$ to $u_i$ with the same index and fixes the
$\Bsort$. Hence the same-code map computes $\Theta$. Since $\Theta$ is an
involution, the same procedure computes its inverse.
\end{proof}

\begin{theorem}\label{thm:infinite-counterexample}
The computable presentation $\Twisted_\infty$ of
$M_\infty^{[0,1)}$ is not aware. Nevertheless, every constant random variable
$\bar n$ is a $\Twisted_\infty$-computable point.
\end{theorem}

\begin{proof}
For each $n\in\N$, consider the rational open ball
\[
 U_n=B\bigl(u_{c(n)};\tfrac13\bigr)
     =B\bigl(\Theta(\bar n);\tfrac13\bigr).
\]
Then $(U_n)_{n\in\N}$ is an effective
sequence of rational open balls in $\Twisted_\infty$.

If $n\in F$, then $q(n)=n$, so $\Theta(\bar n)=\bar n$ and $U_n$ meets the
constant random variables. Suppose that $n\notin F$. Then
$\Theta(\bar n)$ equals $n$ on $I_0$ and $q(n)\neq n$ on $I_1$. Hence, for
every $m\in\N$,
\[
 d_{\Ksort}\bigl(\Theta(\bar n),\bar m\bigr)=
 \begin{cases}
  \tfrac12,&m\in\{n,q(n)\},\\
  1,&m\notin\{n,q(n)\}.
 \end{cases}
\]
Thus $U_n$ is disjoint from the constant random variables. We have proved
\begin{equation}\label{eq:balls-fixed-set}
 U_n\cap\{\bar m:m\in\N\}\neq\varnothing
 \quad\Longleftrightarrow\quad
 n\in F.
\end{equation}

If $\Twisted_\infty$ were aware, then by \Cref{def:ce-closed,def:aware} the
rational open balls meeting $\{\bar m:m\in\N\}$ would form a c.e. set.
Applying this enumeration to the effective sequence $(U_n)$ and using
\Cref{eq:balls-fixed-set} would make $F$ c.e., contradicting
\Cref{lem:infinite-involution}. Therefore $\Twisted_\infty$ is not aware.

It remains to verify the final assertion. For a fixed $n$, the random variable
\[
 v_n(\omega)=\Theta(\bar n)(\omega)=
 \begin{cases}
  n,&\omega\in I_0,\\
  q(n),&\omega\in I_1.
 \end{cases}
\]
is a rational simple random variable and hence occurs among the
$\Ksort$-special points of $\Prob_\infty$. Since $\Theta$ is an involution,
$\Theta(v_n)=\bar n$. Therefore $\bar n$ occurs among the
$\Ksort$-special points of $\Twisted_\infty$ and is
$\Twisted_\infty$-computable. For each fixed $n$, the two natural numbers
$n$ and $q(n)$ can be built into a code for $v_n$; this argument does not
give a uniform procedure for obtaining such a code from $n$.

In fact, there is no uniform procedure producing names
for these constants: otherwise $(\bar n)_{n\in\N}$ would be a
$\Twisted_\infty$-computable enumeration of all constants, contrary to
\Cref{prop:aware-criterion} and the first part of the theorem.
\end{proof}

\section{The proof of Theorem \ref{thm:q36} }

\label{sec:q36}

For convenience we work first over the fair-coin probability space
$(2^\omega,\mu)$.  Binary expansion induces a computable isomorphism between
its probability algebra and the standard probability algebra on $[0,1)$,
modulo the null set of ambiguous binary expansions. 

\begin{definition}[Hirschfeldt \cite{Hirschfeldt}]\label{def:PAC-tree}
\begin{enumerate}
    \item A tree $T\subseteq 2^{<\omega}$ is a \emph{PAC tree} if it is computable, has
no dead ends, and every path through $T$ is computable.
\item If
$S\subseteq2^\omega$, a \emph{computable listing} of $S$ is a uniformly
computable sequence of elements of $2^\omega$ whose range is $S$; repetitions
are allowed.
\item For a PAC tree $T$, let $S_T$ be the set of isolated paths through $T$.
\end{enumerate}
  
\end{definition}

\begin{theorem}[Goncharov--Nurtazin; Millar; see Hirschfeldt {\cite[Theorem~1.9]{Hirschfeldt}}]\label{thm:PAC-tree}
There is a PAC tree $T\subseteq2^{<\omega}$ such that $S_T$ has no computable listing.
\end{theorem}

The next consequence of the PAC property will be used twice. Its proof uses
the path-selection argument from \cite[Theorem~2.1]{Hirschfeldt}.

\begin{lemma}\label{lem:isopath-dense}
Let $T$ be a PAC tree. Then $S_T$ is dense in
$[T]$.
\end{lemma}

\begin{proof}
Fix a noncomputable $X\in2^\omega$ and a node $\sigma\in T$.  Starting at
$\sigma$, follow the unique immediate successor whenever there is only one.
At the $k$th node having two immediate successors, take the successor selected
by $X(k)$.  This defines a path $p\in[T]$ extending $\sigma$.

If infinitely many splitting nodes occurred along $p$, then $X$ would be
computable from $p$: since $T$ is computable, one can locate the successive
splitting nodes along $p$ and read the choices made there.  But every path
through a PAC tree is computable, contradicting the choice of $X$.  Hence
only finitely many splits occur along $p$.  Beyond the last split, $p$ has a
unique continuation in $T$, so $p$ is isolated.
\end{proof}

Fix a tree $T$ as in \Cref{thm:PAC-tree}.  Let
\[
 L_T=\{P_n:n\in\N\}
\]
be a computable language of unary predicates, and define the $L_T$-structure
$$M_T=S_T\times\N,$$
and for all $k$ and $n$, 
\begin{equation}\label{eq:MT-definition}
 M_T\models P_n(p,k)
 \quad\Longleftrightarrow\quad
 p(n)=1.
\end{equation}
$M_T$ is countable since $S_T$ is countable. \Cref{eq:MT-definition} ensures that each path type is
realized infinitely often.

For $\sigma\in2^m$, let
\[
 \vartheta_\sigma(x)=
 \bigwedge_{\substack{j<m\\\sigma(j)=1}}P_j(x)
 \wedge
 \bigwedge_{\substack{j<m\\\sigma(j)=0}}\neg P_j(x).
\]

\begin{lemma}\label{lem:MT-cells}
For every $\sigma\in2^{<\omega}$,
\[
 \vartheta_\sigma(M_T)\neq\varnothing
 \quad\Longleftrightarrow\quad
 \sigma\in T.
\]
Whenever $\vartheta_\sigma(M_T)$ is nonempty, it is infinite.
\end{lemma}

\begin{proof}
If $M_T\models\vartheta_\sigma(p,k)$, then $\sigma\prec p\in[T]$, so
$\sigma\in T$.  Conversely, if $\sigma\in T$, then
\Cref{lem:isopath-dense} gives $p\in S_T$ extending $\sigma$.  Every
$(p,k)$, for $k\in\N$, then satisfies $\vartheta_\sigma$.
\end{proof}

\begin{proposition}\label{prop:MT-QE}
The theory $\operatorname{Th}(M_T)$ is decidable and has effective quantifier
elimination.
\end{proposition}

\begin{proof}
A theory has \emph{quantifier elimination} if every formula is equivalent,
in all models of the theory, to a quantifier-free formula with no additional
free variables; see \cite[Definition~3.1.1]{Marker}.
We use the standard quantifier-elimination argument for a unary structure in
which every nonempty finite Boolean cell is infinite.  To eliminate
$\exists x$ from a quantifier-free formula, first put the formula into
disjunctive normal form.  In any conjunct in which $x$ is required to equal
another variable, substitute that variable for $x$.  Otherwise, $x$ is
required to avoid only finitely many other variables and to satisfy a finite
Boolean pattern in the predicates $P_n$.

Choose $m$ larger than every predicate index appearing in the formula and
consider the finitely many strings $\sigma\in2^m$ extending this partial
pattern.  By \Cref{lem:MT-cells}, the pattern is realizable if and only if at
least one such $\sigma$ belongs to $T$; whenever it is realizable, it has
infinitely many realizations and hence one outside the finite prohibited set.
Membership in $T$ is decidable, so this elimination step is effective.
Iteration gives effective quantifier elimination.  In particular, sentences
can be reduced effectively to quantifier-free sentences and decided.
\end{proof}

\begin{proposition}\label{prop:MT-no-copy}
The structure $M_T$ has no computable presentation.
\end{proposition}

\begin{proof}
Suppose that $(a_i)_{i\in\N}$ were a computable presentation of $M_T$.  For
each $i$, define $p_i\in2^\omega$ by
\[
 p_i(n)=1
 \quad\Longleftrightarrow\quad
 M_T\models P_n(a_i).
\]
The computable atomic diagram makes $(p_i)_{i\in\N}$ a uniformly computable
sequence.  Every $p_i$ lies in $S_T$, and every $p\in S_T$ occurs because the
presentation contains an element corresponding to $(p,0)$.  Thus $(p_i)$ is
a computable listing of $S_T$, contrary to \Cref{thm:PAC-tree}.
\end{proof}

Using the same path-selection procedure, we next construct an
almost-everywhere enumeration of $M_T$ on the probability space $2^\omega$.
Let $(\sigma_i)_{i\in\N}$ be an effective enumeration of
the nodes of $T$ in which every node occurs infinitely often.  Given
$X\in2^\omega$, let $\pi_i^X\in[T]$ be the path obtained by starting at
$\sigma_i$, following the unique successor at a nonsplitting node, and using
the bits of $X$ in order to choose successors at splitting nodes.

\begin{lemma}\label{lem:guided-listing}
If $X$ is noncomputable, then $\pi_i^X\in S_T$ for every $i$, and every
member of $S_T$ occurs infinitely often in $(\pi_i^X)_{i\in\N}$.
\end{lemma}

\begin{proof}
The proof of \Cref{lem:isopath-dense}, applied above $\sigma_i$, shows that
$\pi_i^X$ is isolated: if it passed through infinitely many splitting nodes,
then the choices along the computable path $\pi_i^X$ would compute $X$.

If $p\in S_T$, choose a node $\sigma\prec p$ whose unique extension in
$[T]$ is $p$.  Every guided path starting at $\sigma$ equals $p$.  Since
$\sigma$ occurs infinitely often among the $\sigma_i$, the path $p$ occurs
infinitely often among the $\pi_i^X$.
\end{proof}

For $i\in\N$ and $X\in2^\omega$, define the occurrence rank
\[
 \rho_i(X)=\bigl|\{j<i:\pi_j^X=\pi_i^X\}\bigr|.
\]
Let $\mathsf{Comp}\subseteq2^\omega$ be the set of computable reals.  This set is countable and hence Borel and null.  Fix any $p_*\in S_T$ and define
\begin{equation}\label{eq:fiber-listing}
 f_i(X)=
 \begin{cases}
  (\pi_i^X,\rho_i(X)),&X\notin\mathsf{Comp},\\
  (p_*,i),&X\in\mathsf{Comp}.
 \end{cases}
\end{equation}

\begin{lemma}\label{lem:fiber-bijection}
Each $f_i$ is a Borel $M_T$-valued random variable.  For every
$X\notin\mathsf{Comp}$, the map $i\mapsto f_i(X)$ is a bijection from $\N$
onto $M_T$.
\end{lemma}

\begin{proof}
The map $X\mapsto\pi_i^X$ is continuous: a finite initial segment of its
output is determined by a finite initial segment of $X$.  It follows that
$\{X:\pi_j^X=\pi_i^X\}$ and $\{X:\pi_i^X=p\}$, for fixed $p$, are closed.
Thus $\rho_i$ is Borel, and every fiber of the function in
\eqref{eq:fiber-listing} is Borel. Since $M_T$ is countable and discrete,
this proves that $f_i$ is Borel measurable.

For noncomputable $X$, \Cref{lem:guided-listing} says that every $p\in S_T$
appears infinitely often among the $\pi_i^X$.  Its successive occurrences
have ranks $0,1,2,\ldots$.  Hence every pair $(p,k)\in S_T\times\N$ appears
exactly once among the $f_i(X)$.
\end{proof}

The modification on $\mathsf{Comp}$ makes every $f_i$ genuinely $M_T$-valued
at every point; no effective recognition of $\mathsf{Comp}$ will be needed.
For $i,n\in\N$, put
\[
 E_{i,n}=\{X\in2^\omega:\pi_i^X(n)=1\}.
\]
The set $E_{i,n}$ is effectively clopen, uniformly in $i$ and $n$.  Indeed,
to determine $\pi_i^X(n)$, one extends $\sigma_i$ by at most
$\max\{0,n+1-|\sigma_i|\}$ bits and uses at most one bit of $X$ per
extension. All immediate-successor questions in $T$ are decidable.
Simulating this finite computation on every binary word of that length
therefore produces a finite union of basic cylinders coding $E_{i,n}$.
In the probability algebra we have
\begin{align}
 \angles{P_n(f_i)}&=E_{i,n},
 \label{eq:atomic-predicate}\\
 \angles{f_i=f_j}&=
 \begin{cases}
  2^\omega,&i=j,\\
  \varnothing,&i\neq j.
 \end{cases}
 \label{eq:atomic-equality}
\end{align}
The first equality ignores only the null set $\mathsf{Comp}$.  For the second, the definition on $\mathsf{Comp}$ also keeps
the $f_i(X)$ pairwise distinct there.
Thus equality events are computed by comparing indices; no computation of
the occurrence ranks $\rho_i(X)$ is required.

Let $\mathcal C$ be the effective Boolean algebra of clopen subsets of
$2^\omega$, represented by finite unions of basic cylinders.  Let
$\mathcal D$ be the countable family of all finite clopen pastings of the
$f_i$: its members have the form
\begin{equation}\label{eq:clopen-pasting}
 h(X)=f_{i_a}(X)\quad\text{when }X\in C_a,
\end{equation}
where $(C_a)_{a<s}$ is a finite partition of $2^\omega$ by members of
$\mathcal C$ and $i_0,\ldots,i_{s-1}\in\N$.  The family $\mathcal D$ has an
effective enumeration by the finite codes for these partitions and indices.

\begin{lemma}\label{lem:pasting-dense}
The family $\mathcal D$ is dense in the $\Ksort$-sort of the Borel
randomization of $M_T$ over $2^\omega$.
\end{lemma}

\begin{proof}
Let $g\colon2^\omega\to M_T$ be a Borel random variable.  Outside
$\mathsf{Comp}$, the sequence $(f_i(X))_i$ enumerates $M_T$ without
repetition.  Define pairwise disjoint Borel sets
\[
 A_i=\{X:g(X)=f_i(X)\text{ and }(\forall j<i)\ g(X)\neq f_j(X)\}.
\]
Their union contains $2^\omega\setminus\mathsf{Comp}$.

Given $\eps>0$, choose $N$ such that
$\mu(\bigcup_{i<N}A_i)>1-\eps/2$.  Finite measurable partitions can be
approximated in measure by clopen partitions, so there is a clopen partition
$C_0,\ldots,C_N$ such that
\[
 \sum_{i<N}\mu(A_i\mathbin\triangle C_i)<\eps/2.
\]
Let $h$ equal $f_i$ on $C_i$ for $i<N$ and equal $f_0$ on $C_N$.  Then
\[
 \{X:g(X)\neq h(X)\}
 \subseteq
 \left(2^\omega\setminus\bigcup_{i<N}A_i\right)
 \cup
 \bigcup_{i<N}(A_i\mathbin\triangle C_i),
\]
and hence $d_{\Ksort}(g,h)<\eps$.  Thus $\mathcal D$ is dense. This argument
only establishes density; an effective choice of $N$ or of the approximating
clopen partition is not needed.
\end{proof}

\begin{proposition}\label{prop:MT-randomization-computable}
Taking $\mathcal D$ as the special points of the $\Ksort$-sort and
$\mathcal C$ as the special points of the $\Bsort$-sort yields a computable
presentation of the Borel randomization of $M_T$ over $2^\omega$.
\end{proposition}

\begin{proof}
We show uniformly how to compute all formula events on tuples from
$\mathcal D$.  Given an $L_T$-formula $\varphi(\bar x)$, use
\Cref{prop:MT-QE} to compute a quantifier-free formula $\psi(\bar x)$
equivalent to $\varphi$ in $\operatorname{Th}(M_T)$.  Given a tuple
$\bar h$ of clopen pastings, refine the finitely many partitions defining
its components to a common clopen partition.  On each cell, every component
is equal to some $f_i$.

An equality between two such components is decided by comparing their
indices, using \eqref{eq:atomic-equality}.  An atomic predicate $P_n(f_i)$ is
represented by the uniformly effective clopen set $E_{i,n}$, using
\eqref{eq:atomic-predicate}.  Performing the Boolean operations in $\psi$
therefore produces effectively a clopen set $C_{\varphi,\bar h}$ satisfying
\[
 \angles{\varphi(\bar h)}=C_{\varphi,\bar h}
\]
in the probability algebra.

The measure of a clopen subset of $2^\omega$ is a computable dyadic rational.
Consequently, all formula-event function symbols and their measures are
computable on the named points.  In particular,
\[
 d_{\Ksort}(h,h')=\mu\angles{h\neq h'}
\]
is computable for $h,h'\in\mathcal D$.  Boolean operations, measure, and
$d_{\Bsort}(A,B)=\mu(A\mathbin\triangle B)$ are computable on clopen event
codes, and every event generated from the named points again has an effective
clopen code.  Finally, \Cref{lem:pasting-dense} gives density in the
$\Ksort$-sort, while the clopens are dense in the probability algebra.  Thus
the two families satisfy \Cref{def:computable-presentation}.
\end{proof}

\begin{proof}[Proof of \Cref{thm:q36}]
By \Cref{prop:MT-QE}, the theory of $M_T$ is decidable and has effective
quantifier elimination.  By \Cref{prop:MT-no-copy}, $M_T$ has no computable,
and hence no decidable, presentation.  By
\Cref{prop:MT-randomization-computable} and the computable identification of
the fair-coin probability algebra with the standard probability algebra on
$[0,1)$, the randomization $M_T^{[0,1)}$ has a computable presentation.
\end{proof}

\section*{acknowledgments}
The author was supported by Focused Research Grants of Natural Science Foundation of Jiangsu Province No. BK20243060. The author used OpenAI’s ChatGPT 5.6, to assist with developing the proofs. All proofs in this paper have been independently verified by the author without the use of AI tools. The author takes full responsibility for the mathematical correctness, attribution of ideas, and final content of the paper.

\end{document}